\documentclass{article}
\usepackage{amsmath,amsthm,amssymb,amscd}

\newtheorem{Thm}{Theorem}[section]
\newtheorem{Lem}[Thm]{Lemma}

\theoremstyle{remark}

\theoremstyle{definition}

\def\Mu{\mu}
\def\Mu{\Mu}

\begin{document}

\begin{center}
{\Large \bf Lyapunov exponents of hyperbolic measures and
hyperbolic periodic orbits}\\

\smallskip
Zhenqi Wang \\
\smallskip LMAM, School of Mathematical Sciences, Peking University, Beijing 100871, China \\
\smallskip
E-mail: { \tt wangzq@pku.org.cn}.\\
\smallskip
Wenxiang Sun $^*$ \\
\smallskip
LMAM, School of Mathematical Sciences, Peking University, Beijing 100871, China\\
\smallskip
E-mail: { \tt sunwx@math.pku.edu.cn}.\\
\end{center}

\footnotetext {$^*$ Sun is the corresponding author and is
supported by NNSFC (\# 10231020,
 10671006) and National Basic Research Program of
China(973 Program)(\# 2006CB805900) } 
\footnotetext {AMS Review: 35C50; 37D25}

\begin{abstract}

Lyapunov exponents of a hyperbolic ergodic measure are approximated
by Lyapunov exponents of hyperbolic atomic measures on  periodic
orbits.
\end{abstract}
\bigskip

\section{Introduction}

Let $f$ be a $C^{1+\alpha}$, $\alpha>0,$ diffeomorphism of a
compact $d$-dimensional manifold $M$   and  $Df:\,\,TM \rightarrow
TM$ the derivative of $f.$ Let us fix a smooth Riemannian metric
on $M$, i. e., a scalar product (and consequently a norm) in every
tangent space $T_{x}M$, $x \in M$ which depends on $x$ in a
differentiable way. The limit
    $$\lambda(x, v) =\lim_{n\to \infty}\frac{\log\|Df^{n}v\|}{n},~v\in T_xM,~v \neq
    0,\,\, x\in M\,\,\,\,\,\,\,\,\,\,\,\,\,\,\,\,\,\,\,\,(1.1)$$
is called a Lyapunov exponent for a tangent vector $v\in T_xM$.
Lyapunov exponents describe the asymptotic evolution of tangent
map: positive or negative exponents correspond to exponential
growth or decay of the norm, respectively, whereas vanishing
exponents mean lack of exponential behavior. From the Oseledec
theorem[8] the limit $\lambda(x, v)$ exists for all non vanishing
vectors  $v$ based on almost all  state points $x$ in $M$ with
respect to any given invariant measure, and it is independent of
the  points if the measure is ergodic. The function $\lambda$
being defined on the tangent bundle $TM$ takes on at most $d$
values on each tangent space $T_xM.$ None of these values depend
on the choice of a Riemannian metric. Back to Lyapunov and Perron,
Lyapunov exponents for a differential equation are a natural
generalization of  the eigenvalues of the matrix in  the linear
part of the equation, and the condition that all Lyapunov
exponents are negative together with the Lyapunov-Perron
regularity implies  Lyapunov stability, see Chapter 1 in [1]  for
the definition of Lyapunov-Perron regularity and the classical
theory of Lyapunov stability. The abstract Lyapunov exponent
defined in (1.1) is  a basic concept and an active topic in the
theory of nonuniformly hyperbolic systems known as Pesin theory.
Pesin theory  recovers some hyperbolic behavior for the points
whose Lyapunov exponents are all nonzero. In particular, these
points have well-defined unstable and stable invariant manifolds.
For these reasons, an ergodic invariant measure is called
hyperbolic if its Lyapunov exponents are different from zero.

The closing lemma of Katok[5] in Pesin theory states that
hyperbolic periodic points are dense in the closure of the basin
of a given hyperbolic measure. Based on this lemma we will show in
the present paper that  the Lyapunov exponents of a hyperbolic
ergodic measure are approximated by those   of a hyperbolic atomic
measure on a periodic orbit. Lyapunov exponents for an atomic
measure concentrated on a periodic orbit with period $p$  are
exactly the logarithm  of the norms of eigenvalues for $Df^p.$ Now
we state our main theorem in the present paper.

\begin{Thm}\label{thm:1}
Let $f$  be a  $C^{1+\alpha}$  diffeomorphism of a compact
$d-$dimensional Riemannian manifold $M$ preserving an ergodic
hyperbolic measure $m$ with Lyapunov exponents $\lambda_{1}\leq
\ldots\leq \lambda_{r}<0<\lambda_{r+1}\leq \ldots\leq
\lambda_{d}.$ Then the Lyapunov exponents of $m$ can be
approximated by the Lyapunov exponents of hyperbolic periodic
orbits. More precisely, for any $\gamma>0$, there exists a
hyperbolic periodic point $z$ with Lyapunov exponents
$\lambda_{1}^{z}\leq\ldots\leq\lambda_{d}^{z}$ such that
$\mid\lambda_{i}-\lambda_{i}^{z}\mid<\gamma,$ $ i=1, ..., d.$
\end{Thm}

It is known from  S. Theorem 5.5 in [6]( see also Theorem 15.4.7 in
[2]) that the smallest absolute value in all of the   Lyapunov
exponents of a hyperbolic measure can be approximated from upper
side by the smallest absolute value in all  of the Lyapunov
exponents of a hyperbolic periodic orbit. Our Theorem \ref{thm:1}
could be viewed as a generalization of this  result. Moreover, our
result may contribute to a  strong version of the closing lemma in
[7], which states that a recurrent orbit in a Pesin set of a
hyperbolic measure and the Oseledec splitting it carries can be
approximated by hyperbolic periodic orbits and their Oseledec
splittings. Duo to the discontinuity of the Oseledec splitting, this
is a nontrivial topic.

A classical result of Sigmund[14] in uniform hyperbolic systems
states that periodic measures are dense in the set of invariant
measures. For nonuniform hyperbolic case Hirayama[3] proved that
periodic measures are dense in the set of invariant measures
supported by a total measure set with respect to a hyperbolic
$mixing$ measure. Quite recently, Liang, Liu and Sun[7] replaced the
assumption of hyperbolic $mixing$ measure by a more natural and
weaker assumption of hyperbolic $ergodic$ measure and generalized
Hirayama's result. It is good to point out that this series  of
results of measure approximation do not imply automatically the
present work of exponent approximation.

Using  preliminary  facts recalling in Section 2 we prove in
Section 3 that  the largest(smallest) Lyapuov exponent of an
ergodic hyperbolic measure $m$ can be approximated by that of
atomic measures on hyperbolic periodic orbits. A similar
approximation property  for 2th-exterior power shows that the sum
of the largest(smallest) two Lyapunov exponents of $m$ can be
approximated by that of hyperbolic periodic orbits, which  then
implies   that the two largest(smallest) Lyapunov exponents of $m$
can be approximated by that of hyperbolic periodic orbits.
Inductive  arguments show Theorem \ref{thm:1} in Section 4.

\bigskip

\section{Preliminaries}

In this section, we recall preliminary facts cited  from reference papers and books.

\subsection{A criterion for  hyperbolicity [4]}
We denote by $M$ a $C^\infty$ compact manifold throughout this
paper.  Let $U\subseteq M$ be an open set and let $f$ be a $C^1$
diffeomorphism from $U$ onto $f(U)$. Let $\triangle\subset U$ be a
compact and $f$-invariant set. Let $T_{\Delta}M=E^{1}\oplus E^{2}$
be a Whitney splitting. Put
\\
$$Df=\left[\begin{array}{ccc} G_{11} & G_{12}\\
G_{21} & G_{22}\\
\end{array}\right]:E^{1}\oplus E^{2}\rightarrow E^{1}\oplus
E^{2},$$
$$Df^{-1}=\left[\begin{array}{ccc} G_{11}^{'} & G_{12}^{'}\\
G_{21}^{'} & G_{22}^{'}\\
\end{array}\right]:E^{1}\oplus E^{2}\rightarrow E^{1}\oplus
E^{2},$$ where $G_{ij},\,\,\,\,G'_{ij}$ $i, j=1, 2$ are bundle
maps coving $f.$ If there exist $0<\lambda<1$ and $\epsilon>0$
satisfying $0<\varepsilon<\min\{1-\lambda,\lambda^{-1}-1\}$, and
also\\
$$\max\{\mid G_{11}\mid,\mid G_{22}^{-1}\mid,\mid
G_{11}^{'^{-1}}\mid,\mid G_{22}^{'}\mid\}<\lambda,$$

$$\max\{\mid G_{12}\mid,\mid G_{21}\mid,\mid G_{12}^{'}\mid,\mid
G_{21}^{'}\mid\}<\varepsilon,$$ then $f: \Delta\to \Delta$ is hyperbolic.\\

\subsection{Oseledec theorem[8]}

Let $d$ denote the dimension of the compact manifold $M.$  Let $f:
M\to M$ be a $C^{1}$ diffeomorphism preserving an ergodic
probability measure $m$.
Then there exist\\
(a)  real numbers  $\lambda_{1}< \cdot\cdot\cdot < \lambda_{k}(k \leq d)$;\\
(b)  positive integers $n_{1}, \cdot\cdot\cdot, n_{k}$, satisfying $n_{1}+
\cdot\cdot\cdot +n_{k}=d$; \\
(c)  a Borel set $L(m)$ satisfying $fL(m)=L(m)$ and $m(L(m))=1$;\\
(d)  a measurable splitting
$T_{x}M=E_{x}^{1}\oplus\cdot\cdot\cdot\oplus E_{x}^{k}$
   with dim$E_{x}^{i}=n_{i}$ and $Df(E_{x}^{i})=E_{fx}^{i}$,\\
 such    that
$$\lim_{n\rightarrow\pm \infty}\frac{\log\|Df^{n}v\|}{n}=\lambda_{i},
   $$
   for $\forall x \in L(m)$, $v \in E_{x}^{i},\,\, i=1, 2, \cdots,
   k  $.

The set $L(m)$ is called an  Oseledec basin of $m$.\\

\subsection{Parallelepiped  spectrum (see [13])}
Let $f:M\to M$ be a $C^{1}$ diffeomorphism  preserving an ergodic
measure $m$. The Lyapunov exponents of $m$  being
$\lambda_{1}<\ldots<\lambda_{k}$ with associated splitting
$T_xM=E^{1}\oplus\cdot\cdot\cdot\oplus E^{k},\,\,x\in L(m)$ and
multiplicities $\Gamma(r)=\dim E^{r}$ constitute the spectrum of
$(m,Df)$. We construct a bundle $\wedge^{i}(M),~2\leq i\leq d
$(recall $d=dim M$)  of $C_{d}^{i}$-dimension on $M$, where the
fiber
over $x$ is\\
$$\wedge^{i}(x)=\{v_{j_{1}}\wedge \ldots\wedge v_{j_{i}}:~v_{j_{k}}
\in T_{x}M,\,\,~1\leq k\leq i, ~1\leq j_{1}<j_{2}<\ldots<j_{i}\leq
d\}.$$ Let $Df^{\wedge^{i}}:\Lambda^{i}(M)\rightarrow
\Lambda^{i}(M)$ denote the
$i$-exterior power of $Df$, namely,\\
$$D_{x}f^{\Lambda^{i}}(v_{j_{1}}\wedge \ldots\wedge v_{j_{i}})
=D_{x}f(v_{j_{1}})\wedge\ldots\wedge D_{x}f(v_{j_{i}}).$$ We define
a norm $\parallel~\parallel_{\wedge^{i}}$ on $\wedge^{i}$ by
assigning $v_{j_{1}}\wedge \ldots\wedge v_{j_{i}}$ to $i-$volume of
the parallelepiped generated by the vectors $v_{j_{1}}, ...,
v_{j_{i}}.$ The  spectrum of Lyapunov exponents of  $(m,
Df^{\wedge^{i}})$ consists of numbers
$\varpi=\sum_{r}n_{r}\lambda_{r},$ where $0\leq n_{r}\leq \Gamma(r)$
and $\sum_{r}n_{r}=i$. The subspace corresponding to $\varpi$ in the
associated splitting is generated by $v_{j_{1}}\wedge \ldots\wedge
v_{j_{i}}$, where $v_{j_{l}}\in E^{j_{l}}$ and
$\sum_{l=1}^{i}\lambda^{j_{l}}=\varpi$.

\subsection{Pesin set[9-11]}

 Given $\lambda,\mu \gg \varepsilon>0$, and for all $k \in \mathbb{Z}^{+}$,
we define $\Lambda_{k}=\Lambda_{k}(\lambda,\mu;\varepsilon)$ to be
all points $x \in M$ for which there is a splitting
$T_{x}M=E_{x}^{s} \oplus E_{x}^{u}$ with invariant property
$(D_{x}f^{m})E_{x}^{s}=E_{f^{m}x}^{s}$ and
$(D_{x}f^{m})E_{x}^{u}=E_{f^{m}x}^{u}$ and satisfying:\\
$(a)$ $\|Df^{n}/_{E_{f^{m}x}^{s}}\| \leq e^{\varepsilon
k}e^{-(\lambda-\varepsilon)n}e^{\varepsilon \mid m\mid},~\forall
m\in\mathbb{Z},~n\geq 1$;\\
$(b)$ $\|Df^{-n}/_{E_{f^{m}x}^{u}}\| \leq e^{\varepsilon
k}e^{-(\mu-\varepsilon)n}e^{\varepsilon \mid m\mid},~\forall m\in
\mathbb{Z}, ~n\geq 1$;\\
$(c)$ $\tan (Angle(E_{f^{m}x}^{s},E_{f^{m}x}^{u})) \geq
e^{-\varepsilon k}e^{-\varepsilon \mid m\mid},~\forall m\in
\mathbb{Z}$. \\
We put
$\Lambda=\Lambda(\lambda,\mu;\varepsilon)=\bigcup_{k=1}^{+\infty}
\Lambda_{k}$ and call $\Lambda$ a Pesin set.

Let $m$ be an   ergodic hyperbolic measure preserved by $f$. We
denote by $\lambda$ the absolute value of the largest negative
Lyapunov exponent and $\mu$ the smallest positive Lyapunov
exponent of $m.$ Let $E^s$ and $E^u$ denote, respectively, the
direct sum of the sub bundles corresponding to negative  Lyapunov
exponents and the sum of sub bundles corresponding to positive
exponents. Then $E^s$ and $E^u$ are well defined on the Oseledec
basin L(m)( see 2.2), they are $Df$ invariant and their direct sum
based on L(m) coincides with $T_{L(m)}M.$ By using these $\lambda$
and $\mu$ together with $E^s$ and $E^u$ we get as in the above
definition a Pesin set $\Lambda=\Lambda(\lambda,\mu;\varepsilon)$
for small $\varepsilon$. This is called a Pesin set of $m.$
 It follows that $m(\Lambda\setminus L(m))+
m(L(m)\setminus \Lambda)=0.$

The following statements are elementary: \\
(1) $\Lambda_{1} \subseteq \Lambda_{2} \subseteq \Lambda_{3}
\cdot\cdot\cdot\cdot$;\\
(2) $f(\Lambda_{k}) \subseteq \Lambda_{k+1},~f^{-1}(\Lambda_{k})
\subseteq \Lambda_{k+1}$;\\
(3) $\Lambda_{k}$ is compact for $\forall k\geq 1$;\\
(4) for $\forall k\geq 1$ £¬the splitting $x\to E_{x}^{u}\oplus
E_{x}^{s}$ depends continuously on  $\Lambda_{k}$.\\

\subsection{Lyapunov metric $\|~\|^{'}$ [9-11]}

Let $\lambda^{'}=\lambda-2\varepsilon,~ \mu{'}=\mu-2\varepsilon.$
Note that $\epsilon<< \lambda, \mu,$ then $\lambda', \mu'>0.$ Let
$x\in \Lambda(\lambda, \mu, \epsilon),$ a Pesin set, see 2.4. For
$v_{s}\in E_{x}^{s},$ we define
$\|v_{s}\|_{s}=\sum_{n=0}^{+\infty}e^{\lambda^{'}n}\|D_{x}f^{n}(v_{s})\|$;
for $v_{u}\in E_{x}^{u},$ we define
$\|v_{u}\|_{u}=\sum_{n=0}^{+\infty}e^{\mu
^{'}n}\|D_{x}f^{-n}(v_{u})\|$. And we define Lyapunov metric $\|
\,\,\|^{'}$ on $T_\Lambda M$ by $\|v\|^{'}=\max\{\|v_{s}\|_{s},
\|v_u\|_{u}\},$ where  $v=v_{s}+v_{u} \in E_{x}^{s}\oplus
E_{x}^u$, $x\in \Lambda.$ We call as usual the norm $\| .\|' $ a
Lyapunov metric. This metric is in general not equivalent to the
Riemannain metric. With the Lyapunov metric  $f:\Lambda
\rightarrow \Lambda$ is uniformly hyperbolic.  The following estimates  are  known :\\
$(a)$ $\|Df/_{E_{x}^{s}}\|^{'} \leq e^{-\lambda^{'}},
~\|Df^{-1}/_{E_{x}^{u}}\|^{'} \leq e^{-\mu^{'}}$;\\
$(b)$ $\frac {1}{\sqrt{d}}\|v\|_{x}\leq \|v\|_{x}^{'} \leq
Ce^{\varepsilon k}\|v\|_{x}, ~\forall v\in T_{x}M, ~x\in
\Lambda_{k}$, where $C=\frac {2}{1-e^{-\varepsilon}}$.\\

In 2.1 through 2.5, the diffeomorphism $f$ is  supposed to be of
$C^1.$ From now on, we assume that $f$ is $C^{1+\alpha},$
$0<\alpha<1,$ that is, $f$ is $C^1$ and furthermore there exists a
costant $K>0$ so that
$$\parallel D_xf-D_yf\parallel \leq K d(x, y)^\alpha,\,\,\,\, \forall x,\,y\in M,$$
provided  $d(x, y)$ is small.

\subsection{Extension of Lyapunov metric[9-11]}

 Fix a point $x\in \Lambda=\Lambda(\lambda, \mu, \epsilon),$
where $\Lambda(\lambda, \mu, \epsilon)$ is a Pesin set, see 2.4.
By taking charts about $x,~f(x)\in M$ we can assume without loss
of generality that $x\in \mathbb{R}^{d},~f(x)\in \mathbb{R}^{d}$.
For a sufficiently small neighborhood $U$ of $x$, we can
trivialize the tangent bundle over $U$ by identifying
$T_{U}M\equiv U\times \mathbb{R}^{d}$. For any point $y\in U$ and
tangent vector $v\in T_{y}M$ we can then use the identification
$T_{U}M\equiv U\times \mathbb{R}^{d}$ to $translate$ the vector
$v$ to a corresponding vector $\overline{v}\in T_{x}M$. We then
define $\parallel v\parallel_{y}^{''}=\parallel
\overline{v}\parallel_{x}^{'},$ where $\|.\|'$ indicates  the
Lyapunov metric defined in 2.5. This defines a new norm
$\|~\|^{''}$ on $T_{U}M$ (which agrees with $\|~\|^{'}$ on the
fiber $T_{x}M$). Similarly, we can define $\|~\|_{z}^{''}$ on
$T_{z}M$ (for any $z$ in a sufficiently small neighborhood  of
$fx$). We write $\overline{v}$ as $v$ whenever there is no
confusion. We can define a new splitting
$T_{y}M=E_{y}^{s^{'}}\oplus E_{y}^{u^{'}},$  $y\in U$ by
$translating$ the splitting $T_{x}M=E_{x}^{s}\oplus E_{x}^{u}$
(and similarly for $T_{z}M=E_{z}^{s^{'}}\oplus E_{z}^{u^{'}}$).

There exist $0<\lambda^{''}<\lambda^{'},~0<\mu^{''}<\mu^{'}$ and
$\varepsilon_{0}>0$ such that if we set
$\varepsilon_{k}=\varepsilon_{0}e^{-\varepsilon k}$ then for any
$y\in B(x,\varepsilon_{k})$ in an $\varepsilon_{k}$ neighborhood
of $x\in \Lambda_{k}$ we have a splitting
$T_{y}M=E_{y}^{s^{'}}\oplus E_{y}^{u^{'}}$ with hyperbolicity behavior: \\
$\|D_{y}f(v)\|_{fy}^{''} \leq e^{-\lambda^{''}}\|v\|_{y}^{''}$,
for every $v\in
E_{y}^{s^{'}};$\\
$\|D_{y}f^{-1}(w)\|_{f^{-1}y}^{''}\leq e^{-\mu^{''}}
\|w\|_{y}^{''}$, for every $w\in
E_{y}^{u^{'}}.$\\
The constant $\epsilon_0$ here and afterwards depends on various
global properties of $f,$ e. g. the Holder constants, the size of
the local trivialization, see p.73 in [12].

\subsection{Shadowing lemma and closing lemma}
Let $(\delta_{k})_{k=1}^{+ \infty}$ be a sequence of positive real
numbers. Let $(x_{n})_{n=-\infty}^{+ \infty}$ be a sequence in
$\Lambda=\Lambda(\lambda, \mu, \epsilon)$ for which there exists a
sequence $(s_{n})_{n=-\infty}^{+ \infty}$ of positive integers
satisfying:\\
$(a)$ $x_{n}\in \Lambda _{s_{n}}, ~\forall n\in \mathbb{Z}$;\\
$(b)$ $\mid s_{n}-s_{n-1}\mid \leq 1, ~\forall n\in \mathbb{Z}$;\\
$(c)$ $d(fx_{n},x_{n+1})\leq \delta_{s_{n}}, ~\forall n\in \mathbb{Z}$;\\
then we call $(x_{n})_{n=-\infty}^{+ \infty}$ a
$(\delta_{k})_{k=1}^{+ \infty}$  pseudo-orbit. Given $\eta>0$, a
point $x\in M$ is an
 $\eta$-shadowing point for the $(\delta_{k})_{k=1}^{+ \infty}$
 pseudo-orbit if $d(f^{n}x,x_{n})\leq \eta \varepsilon_{s_{n}},~ \forall n\in
\mathbb{Z}$, where $\varepsilon_{k}=\varepsilon_{0}e^{-\varepsilon
k}.$  \\

\begin{Lem}
(Shadowing lemma [5] [12, Thm. 5.1]) Let $f:M\rightarrow M$ be a
$C^{1+\alpha}$ diffeomorphism, with a non-empty Pesin set
$\Lambda=\Lambda(\lambda,\mu;\varepsilon)$ and fixed parameters,
$\lambda,\mu\gg\varepsilon>0$. For $\forall \eta >0$ there exists a
sequence $(\delta_{k})_{k=1}^{+ \infty}$ such that for any
$(\delta_{k})_{k=1}^{+ \infty}$ pseudo-orbit there exists a unique
$\eta$-shadowing
point. \\
\end{Lem}

{\bf Remark} If we change $\varepsilon_{0}e^{-ak}$ for
$\varepsilon_{k}=\epsilon_0 e^{-\epsilon k},$  where
 $\min(\lambda-2\epsilon,\, \mu-2\varepsilon)>a\geq \varepsilon$,
 then the shadowing lemma is still
 true(see the argument at  pages. 89-93 in [12]).\\

\begin{Lem}
(Closing lemma[5]) Let $f:M\rightarrow M$   be a $C^{1+\alpha}$
diffeomorphism and let $\Lambda=\Lambda(\lambda,\mu;\varepsilon)$ be
a non-empty Pesin set. For $\forall k \geq 1,\,\, 0<\eta<1$, there
exists $\beta=\beta(k,\eta)>0$ such that: if $x,f^{p}x \in
\Lambda_{k}$ and $d(x,f^{p}x)<\beta$ then there exists a periodic
point
$z\in M,$ with $z=f^{p}z$ and $d(z,x)<\eta$.\\
\end{Lem}

{\bf Remark} By shadowing lemma(Lemma 2.1) and its remark, we
easily get a more convenient version of the
closing lemma as follows:\\
 Let $f:M\rightarrow M$  be a $C^{1+\alpha}$  diffeomorphism and let
 $\Lambda=\Lambda(\lambda,\mu;\varepsilon)$  be a
non-empty Pesin set. For $\forall k \geq 1, 0<\eta<1$,
$\min(\lambda-2\varepsilon,~\mu-2\varepsilon)> \theta\geq
\varepsilon$, there exists $\beta=\beta(k,\eta,\theta)>0$ with the
property that  if $x,f^{p}x \in \Lambda_{k}$ and
$d(x,f^{p}x)<\beta$ then there exists a periodic point $z\in M,$
$z=f^{p}z$, such that
$d(f^{i}x,f^{i}z)<\eta \varepsilon_{0}e^{-\theta i}$, for
$0\leq i \leq p-1$ (see p. 95 in [12]).\\

\section {The largest and the  smallest Lyapunov exponents}

In this section we show that  the largest and the smallest Lyapunov
exponents of an ergodic  hyperbolic measure are approximated by that
of hyperbolic periodic orbits.

\begin{Thm}\label{thm:2}
 Let
$f:M \rightarrow M$ be a $C^{1+\alpha},$ $0<\alpha<1,$
diffeomorphism of a compact manifold of dimension $d,$ and let $m$
be an ergodic hyperbolic measure with Lyapunov exponents
$\lambda_{1}<\cdots<\lambda_{r}<0<\lambda_{r+1}<\cdots<\lambda_{t}(t
\leq d)$. Then the largest  Lyapunov exponent $\lambda_t$ can be
approximated by the largest  Lyapunov exponents of hyperbolic
periodic orbits. More precisely, for any $\gamma>0$, there exists a
hyperbolic periodic point $z$ with Lyapunov exponents
$\lambda_{1}^{z}\leq\ldots\leq \lambda_{d}^{z}$ such that
$\mid\lambda_{t}-\lambda_{d}^{z}\mid<\gamma$.
\end{Thm}

Before proving  Theorem \ref{thm:2} we  explain  the main  idea. By
using Katok's shadowing lemma  we get a hyperbolic periodic orbit to
trace a certain  segment of orbit in a Pesin set of $m.$ The orbit
segment is uniformly hyperbolic under Lyapunov metric and the norm
of $Df$(which relates closed to the largest exponent) when
restricted on the orbit segment can be controlled very well by the
Lyapunov exponent under suitable Lyapunov metric. Now that the
periodic orbit is in a small neighborhood of the  orbit segment it
traces and $f$ is $C^{1+\alpha}$, we then transfer the counting
Lyapunove exponent from  the periodic orbit to the orbit segment.
This enables us  to compare the two largest Lyapunov exponents and
to estimate their difference.

The proof is somehow technical. We define one  Pesin set for  $m$
and three Lyapunov metrics for  the Pesin set  by using all the
individual  exponents and thus all the corresponding individual
sub bundles in the Oseledec splitting, comparing with the standard
Pesin set as  in 2.4 and the standard  Lyapunov metric as in 2.5
by using the largest negative exponent and the smallest positive
exponent and thus the stable bundle that is the direct sum of sub
bundles corresponding to all  negative exponents and the unstable
bundle that is the direct sum of the sub bundles corresponding to
all positive exponents. One  of the advantages of our definitions
is that they enable us to control  the norm of derivative
restricted on each sub bundle by corresponding  exponent from both
lower side and upper side. Another advantage is  that we get three
pairs of desired estimates (3.1.1)-(3.1.2), (3.2.1)-(3.2.2) and
(3.3.1)-(3.3.2) under  new metrics, comparing with the
inequalities (a) (b) under standard Lyapunov metric as in 2.5. By
using (3.1.1)-(3.1.2) and Katok's closing lemma and criterion in
2.1  we prove  the existence of a periodic orbit $orb(z)$ which is
hyperbolic under the first Lyapunov metric we defined in the
proof. (3.2.1)-(3.2.2) contribute to proving that
$$\lim_{n\rightarrow +\infty}\frac {\log\parallel
D_{z}f^{n}\parallel^{(4)}}{n}<\lambda_{t}+\gamma,$$ where
$\parallel \cdot \parallel^{(4)}$ denotes the extension metric to
the second  Lyapunov metric we defined in the proof.
(3.2.3)-(3.3.2) contribute to proving that
$$\lambda_t-\gamma< \lim_{n\rightarrow +\infty}\frac
{\log\parallel D_{z}f^{n}\parallel^{(6)}}{n},$$ where $\parallel
\cdot \parallel^{(6)}$ denotes the extension metric to  the third
Lyapunov metric we defined in the proof. These inequalities give
rise to the final inequality
$$\mid\lambda_{t}-\lambda_{d}^{z}\mid<\gamma$$
under the Riemannian metric, although what we defined three metrics
are not equivalent to the Riemannian one on a whole Pesin set in
general. This is because in our case only a finite number of Pesin
Blocks are used and thus the Lyapunov metrics restricted on these
blocks are equivalent to the Riemannian one.

The proof is divided into three steps. The existence of hyperbolic
orbit $Orb(z)$  in step 1 is not a new result, which was  proved
in Katok[5]. And our   proof is not quite different from in [5].
But it is better adapted to the proof of the following two steps.

\bigskip
{\bf Proof of Theorem \ref{thm:2} } Given $\min_{1\leq i\not=j\leq
t}\mid \lambda_{i}-\lambda_{j}\mid \gg \varepsilon>0$, and for all
$k \in \mathbb{Z}^{+}$, we define
$\Lambda_{k}=\Lambda_{k}(\{\lambda_{1},\ldots,\lambda_{t}\};\varepsilon)
$ to be all points $x \in M$ for which there is a splitting
$T_{x}M=E_{x}^{1}\oplus \cdot\cdot\cdot\oplus E_{x}^{t}$ with
$$\lim_{n\to \infty}\frac{\log\|Df^{n}v\|}{n}=\lambda_{i},~0 \neq
v\in E_{x}^{i}$$ and  with invariant property
$(D_{x}f^{m})E_{x}^{i}=E_{f^{m}x}^{i},~1\leq i\leq t$
and satisfying:\\
\indent$(a)$ $e^{-\varepsilon
k}e^{(\lambda_{i}-\varepsilon)n}e^{-\varepsilon \mid
m\mid}\leq\|Df^{n}/_{E_{f^{m}x}^{i}}\| \leq e^{\varepsilon
k}e^{(\lambda_{i}+\varepsilon)n}e^{\varepsilon \mid m\mid},~1\leq
i\leq t,~\forall m\in
\mathbb{Z},~n\geq 1$;\\

$(b)$ $\tan (Angle(E_{f^{m}x}^{i},E_{f^{m}x}^{j})) \geq
e^{-\varepsilon k}e^{-\varepsilon \mid m\mid},~\forall i \neq
j,~\forall m\in
\mathbb{Z}$. \\

We set
$\Lambda=\Lambda(\{\lambda_{1},\ldots,\lambda_{t}\};\varepsilon)=\bigcup_{k=1}^{+\infty}
\Lambda_{k}$ and call $\Lambda$ a Pesin set. We easily get that
$m(\Lambda)=1.$ This  Pesin set is slightly different from the
standard one as in 2.4, but the properties(1)-(4) stated there are
still true.

Let $q=~\min_{1\leq i\not=j\leq t}\mid
\lambda_{i}-\lambda_{j}\mid$, and take arbitrarily $\gamma>0$,
satisfying $\min\{\frac 12,~\frac {5}{d},\frac {q}{2},
\lambda_t\}>\gamma>0$, and satisfying $\log\frac {\frac
{2(t-1)\gamma}{5}+1}{(1-\frac {2\gamma}{5})(1-\frac
{t\gamma}{5})}<q$.
  Let $\varepsilon\leq
\min\{\frac {\gamma}{5},~\frac {1}{4}q-\frac {1}{4}\log\frac
{\frac {2(t-1)\gamma}{5}+1}{(1-\frac {2\gamma}{5})(1-\frac
{t\gamma}{5})}\}$ and $q\gg \varepsilon>0$.

We divide the proof into three steps.
\bigskip

{\bf Step 1} We prove the existence of hyperbolic periodic points
 near our  Pesin set   $\Lambda.$ Although the way  that Katok[5] proved  the
 existence of hyperbolic periodic points near the standard Pesin
 set as in 2.4 works in our case here, we  present a short proof of the
 existence  of periodic
 points near
 our Pesin set by  shadowing lemma and prove the hyperbolicity of these
 periodic points by  the criterion 2.1, a slight different method from that in [5]. The
 techniques and
 inequalities developed while proving the hyperbolicity turns out to be
 helpful to  our consecutive steps.

Let  $\lambda_{i}^{'}=\mid\lambda_{i}\mid-2\varepsilon.$ Then
$\lambda_{i}^{'}>0.$ We define  a new norm $\|v_{i}\|_{i}$ on the
spaces $E_{x}^{i},$ $1\leq i\leq t,$ $x\in \Lambda.$ For $v_{i}\in
E_{x}^{i}, ~1\leq i\leq r$,  we define
$\|v_{i}\|_{i}=\sum_{n=0}^{+\infty}e^{\lambda_{i}^{'}n}\|D_{x}f^{n}(v_{i})\|$;
for $v_{j}\in E_{x}^{j},~r+1\leq j\leq t$, we define
$\|v_{j}\|_{j}=\sum_{n=0}^{+\infty}e^{\lambda_{j}^{'}n}\|D_{x}f^{-n}(v_{j})\|$.
All these series are convergent. For
$v=\sum_{i=1}^{t}v_{i},~v_{i}\in E_{x}^{i}$, we define
$\|v\|^{(1)}=\max_{1\leq i\leq t}\{\|v_{i}\|_{i}\}$. The norm
$\|\,\, \|^{(1)}$  is called   in the present paper Lyapunov metric
number 1, which coincides with the standard Lyapunov metric in 2.5
when $dim M\leq 2.$ This is not equivalent to the Riemannian metric
in general.  With this norm, $f:\Lambda
\rightarrow \Lambda$ is uniformly hyperbolic. The following estimates are similar to that in 2.5. \\
$\|Df/_{E_{x}^{i}}\|^{(1)} \leq e^{-\lambda_{i}^{'}},~1\leq i\leq
r,
~\|Df^{-1}/_{E_{x}^{j}}\|^{(1)} \leq e^{-\lambda_{j}^{'}},~r+1\leq j\leq t;\quad\cdots(3.1.1)$\\
$\frac {1}{\sqrt{d}}\|v\|_{x}\leq \|v\|_{x}^{(1)} \leq
Ce^{\varepsilon k}\|v\|_{x}, ~\forall v\in T_{x}M, ~x\in
\Lambda_{k},\quad\cdots(3.1.2)$\\
 where $C=\frac
{2}{1-e^{-\varepsilon}}$. By  2.6 one extends this norm  to a norm
$\parallel~\parallel^{(2)}.$

From continuity of the Riemannian metric, there exists $\delta>0$
such that
$$\frac 1{1+\gamma}<\frac {\|~\|_{x}}{\|~\|_{y}}<1+\gamma,\quad \cdots(3.1.3)$$
provided $d(x,y)<\delta.$
 Fix $\alpha_{0}$ with
$\min(\lambda_{r+1}-2\varepsilon,\mid\lambda_{r}\mid-2\varepsilon)>\alpha_{0}>2\frac
{\varepsilon}{\alpha}$, where $\alpha$ is the H$\ddot{o}$lder
constant of $Df$. Since
$m(\Lambda(\{\lambda_{1},\ldots,\lambda_{t}\};\varepsilon))=1$,
there exists $k_{0}\in \mathbb{N}$ such that $m(\Lambda_{k_{0}})>0$.
For a given arbitrarily $\eta>0$  we choose
$\beta=\beta(k_{0},\eta,\alpha_{0})>0$ as in the remark to Lemma 2.2
with $\beta < \delta.$ There exists $y_{0}\in \Lambda_{k_{0}}$ such
that $m(B(y_{0},\frac {\beta}{2})\bigcap \Lambda_{k_{0}})>0$ by
compactness of $\Lambda_{k_{0}}$. By Poincar$\acute{e}$'s recurrence
theorem, $\exists \,\, y\in B(y_{0},\frac {\beta}{2})\bigcap
\Lambda_{k_{0}},$  $\exists \,\, p>1$ such that $f^{p}y\in
B(y_{0},~\frac {\beta}{2})\bigcap\Lambda_{k_{0}}.$ Since
$d(y,f^{p}y)<\beta$ and $\min(\lambda_{r+1}-2\epsilon,
|\lambda_r|-2\epsilon )>\alpha_{0}>2\frac
{\varepsilon}{\alpha}>\varepsilon$, by Lemma 2.2 and its remark
there exists a periodic point $z\in M,~z=f^{p}z$,
with\\
 $$d(f^{i}y,f^{i}z)<\eta \varepsilon_{0}e^{-\alpha_{0}i}, \,\, ~0\leq i \leq
 p-1.\quad\cdots\cdots (3.1.4)$$
For $\forall v\in T_{f^{i}y}M,~0\leq i\leq p-2$, by definition of $\|~\|^{(2)}$ we have\\
$$\parallel D_{f^{i}z}fv-D_{f^{i}y}fv
\parallel_{f^{i+1}y}^{(2)}=\parallel D_{f^{i}z}fv-D_{f^{i}y}fv
\parallel_{f^{i+1}y}^{(1)}.$$
Using (3.1.2), (3.1.4) and noting $f^{i+1}y\in \Lambda_{k_{0}+i+1}$,
 we have\\
\begin{eqnarray*}
&\parallel D_{f^{i}z}fv-D_{f^{i}y}fv
\parallel_{f^{i+1}y}^{(1)}\\
\leq &Ce^{\varepsilon(k_{0}+i+1)}\| D_{f^{i}z}fv-D_{f^{i}y}fv
\|_{f^{i+1}y}\\
 \leq &Ce^{\varepsilon(k_{0}+i+1)}Kd(f^{i}z,f^{i}y)^{\alpha}\|v\|_{f^{i}y}\\
\leq&C\sqrt{d}e^{\varepsilon(k_{0}+1)}K\eta^{\alpha}\varepsilon_{0}^{\alpha}
e^{-(\alpha_{0}\alpha-\varepsilon)i}\|v\|_{f^{i}y}^{(2)},\\
\end{eqnarray*}
and thus
$$\parallel D_{f^{i}z}f-D_{f^{i}y}f
\parallel^{(2)}
\leq
C\sqrt{d}e^{\varepsilon(k_{0}+1)}K\eta^{\alpha}\varepsilon_{0}^{\alpha}
e^{-(\alpha_{0}\alpha-\varepsilon)i},~0\leq i \leq
p-2.\cdots(3.1.5)$$

Similarly, for $f^{-1}$ and  $v\in T_{f^{i}y}M,$  $1\leq i\leq p-1,$
we have by (3.1.2) and (3.1.4)\\
\begin{eqnarray*}
&\parallel D_{f^{i}z}f^{-1}v-D_{f^{i}y}f^{-1}v
\parallel_{f^{i-1}y}^{(2)}\\
\leq & Ce^{\varepsilon(k_{0}+i-1)}\parallel
(D_{f^{i-1}z}f)^{-1}(D_{f^{i-1}y}f-D_{f^{i-1}z}f)(D_{f^{i-1}y}f)^{-1}(v)
\parallel_{f^{i-1}y}\\
\leq& Ce^{\varepsilon(k_{0}+i-1)}K\sqrt{d}\parallel
(D_{f^{i-1}z}f)^{-1} \parallel\parallel (D_{f^{i-1}y}f)^{-1}
\parallel d(f^{i-1}z,f^{i-1}y)^{\alpha}\parallel v
\parallel_{f^{i}y}^{(1)}\\
\leq &Ce^{\varepsilon k_{0}}K\sqrt{d}\eta^{\alpha}\| Df^{-1}
\|^{2}\varepsilon_{0}^{\alpha}e^{-(\alpha_{0}\alpha-\varepsilon)(i-1)}\parallel
v
\parallel_{f^{i}y}^{(2)},\\
\end{eqnarray*}
and thus
$$\parallel D_{f^{i}z}f^{-1}-D_{f^{i}y}f^{-1}
\parallel^{(2)}\leq Ce^{\varepsilon k_{0}}K\sqrt{d}\eta^{\alpha}\| Df^{-1}
\|^{2}\varepsilon_{0}^{\alpha}e^{-(\alpha_{0}\alpha-\varepsilon)(i-1)},~1\leq
i\leq p-1.\cdots(3.1.6)$$

For  $v\in T_{f^{p-1}y}M$ by (3.1.3) we have\\
\begin{eqnarray*}
&\parallel D_{f^{p-1}z}fv-D_{f^{p-1}y}fv
\parallel_{f^{p}y}^{(2)}\\
\leq &Ce^{\varepsilon k_{0}}\parallel
D_{f^{p-1}z}fv-D_{f^{p-1}y}fv
\parallel_{y}\\
\leq &(1+\gamma)Ce^{\varepsilon k_{0}}\parallel
D_{f^{p-1}z}fv-D_{f^{p-1}y}fv
\parallel_{f^{p}y}\\
\leq &(1+\gamma)KCe^{\varepsilon
k_{0}}d(f^{p-1}z,f^{p-1}y)^{\alpha}
\parallel v\parallel_{f^{p-1}y}\\
\leq
&(1+\gamma)\sqrt{d}KC\eta^{\alpha}\varepsilon_{0}^{\alpha}e^{\varepsilon
k_{0}}e^{-\alpha_{0}\alpha(p-1)}
\parallel v\parallel_{f^{p-1}y}^{(2)},\\
\end{eqnarray*}
and thus
$$\parallel D_{f^{p-1}z}f-D_{f^{p-1}y}f
\parallel^{(2)}\leq (1+\gamma)\sqrt{d}KC\eta^{\alpha}\varepsilon_{0}^{\alpha}e^{\varepsilon
k_{0}}e^{-\alpha_{0}\alpha(p-1)}.\cdots(3.1.7)$$

For $v\in T_{f^{p}y}M$ we have by using similar estimates as above\\
$$\parallel D_{z}f^{-1}-D_{f^{p}y}f^{-1}
\parallel^{(2)}\leq (1+\gamma)C\eta^{\alpha}\varepsilon_{0}^{\alpha}e^{\varepsilon
k_{0}}K\sqrt{d}\| Df^{-1}
\|^{2}e^{-(\alpha_{0}\alpha-\varepsilon)(p-1)}.\cdots(3.1.8)$$

Let
\begin{eqnarray*}
B=\max\{ &Ce^{\varepsilon(k_{0})}K\sqrt{d}\| Df^{-1}
\|^{2}\varepsilon_{0}^{\alpha}, \,\,
(1+\gamma)C\varepsilon_{0}^{\alpha}e^{\varepsilon
k_{0}}K\sqrt{d}\| Df^{-1} \|^{2},\\
 &C\sqrt{d}e^{\varepsilon(k_{0}+1)}K\varepsilon_{0}^{\alpha},\,\,
 (1+\gamma)
\sqrt{d}KC\varepsilon_{0}^{\alpha}e^{\varepsilon k_{0}}\}.
\,\,\,\,\,\, \cdots\cdots\cdots\cdots\cdots(3.1.9)
\end{eqnarray*}

From (3.1.5)-(3.1.9) we have
$$\parallel
D_{f^{i}z}f-D_{f^{i}y}f\parallel^{(2)}\leq
B\eta^{\alpha}e^{-(\alpha_{0}\alpha-\varepsilon)i}, ~0\leq i \leq
p-1,\cdots (3.1.10)$$ and
$$\parallel
D_{f^{i}z}f^{-1}-D_{f^{i}y}f^{-1}
\parallel^{(2)}\leq
B\eta^{\alpha}e^{-(\alpha_{0}\alpha-\varepsilon)(i-1)},~ 1\leq i
\leq p,\cdots (3.1.11)$$ where we remember from the choice of
$\alpha_{0}$ that $\alpha_{0}\alpha-\varepsilon>0.$

Let  $$E_{x}^{s}=E_{x}^{1}\oplus\cdots\oplus
E_{x}^{r},\,\,\,\,~E_{x}^{u}=E_{x}^{r+1}\oplus\cdots\oplus
E_{x}^{t}.\,\,\cdots \cdots (3.1.12)$$ Let
$A=\max(e^{\lambda_{r}+2\varepsilon},e^{-\lambda_{r+1}+2\varepsilon})$.
Consider  a system of inequations:\\
$$\left\{\begin{array}{ccc} x<1-A-x\\
x<\frac {1}{A+x}-1 \end{array}\right.$$
or
$$\left\{\begin{array}{ccc} x<\frac {1-A}{2}\\
x^2+(1+A)x+A-1<0.\\
\end{array}\right.$$ Since  $A<1$,  there exists a real
number $b>0$ such that any number included in $(0,b)$ is a
solution of the system.  Now we make a  restriction that  $\eta\leq (\frac
{b}{B})^{\frac {1}{\alpha}}$( we will make another restriction in Step 2),
then $B\eta^{\alpha}$ is a solution of the system, i.e.\\
$$\left\{\begin{array}{cc} B\eta^{\alpha}<1-(A+B\eta^{\alpha})\\
B\eta^{\alpha}<\frac {1}{A+B\eta^{\alpha}}-1\\
\end{array}
\right. \cdots (3.1.13)$$ Under the invariant splitting
$E_{f^{i}y}^{s} \oplus E_{f^{i}y}^{u}~i\in \mathbb{Z}$,
$D_{f^{i}y}f$ and $D_{f^{i}y}f^{-1}$  are diagonal block matrixes,
and
$D_{f^{i}z}f$ and $D_{f^{i}z}f^{-1}$ are block matrixes as follows\\
$$D_{f^{i}z}f=\left[\begin{array}{ccc}  {G_{11}^{i}} &
{G_{12}^{i}}\\
{G_{21}^{i}} &
{G_{22}^{i}}\\
\end{array}
\right ]: E_{f^{i}y}^{s} \oplus E_{f^{i}y}^{u} \rightarrow
E_{f^{i+1}y}^{s} \oplus E_{f^{i+1}y}^{u},$$
$$D_{f^{i}z}f^{-1}=\left[\begin{array}{ccc}  {G_{11}^{i^{'}}} &
{G_{12}^{i^{'}}}\\
{G_{21}^{i^{'}}} &
{G_{22}^{i^{'}}}\\
\end{array}
\right ]: E_{f^{i+1}y}^{s} \oplus E_{f^{i+1}y}^{u} \rightarrow
E_{f^{i}y}^{s} \oplus E_{f^{i}y}^{u}.$$ By
(3.1.1),\,\,(3.1.10)-(3.1.12), we have
$$\max(\parallel G_{11}^{i}\parallel^{(2)},\parallel G_{22}^{i^{-1}} \parallel ^{(2)},
 \parallel G_{11}^{i^{'^{-1}}}\parallel^{(2)},
\parallel G_{22}^{i^{'}} \parallel^{(2)} )<A+B\eta^{\alpha},$$
$$\max(\parallel G_{21}^{i}\parallel^{(2)} , \parallel G_{21}^{i}\parallel^{(2)}
 , \parallel G_{12}^{i^{'}}\parallel^{(2)},
\parallel G_{21}^{i^{'}}\parallel^{(2)} )<B\eta^{\alpha}.$$
 According to (3.1.13) and the criterion 2.1,
$orb(z,f)$ is uniformly hyperbolic with  the norm $\|~\|^{(2)}$.
Observe that  $orb(z,f)$ consists of  finite many  points, $orb(z,
f)$ is hyperbolic as well with the Riemannian  norm $\|~\|$.

\bigskip

{\bf Step 2} We prove that $\lim_{n\rightarrow +\infty}\frac
{\log\parallel D_{z}f^{n}\parallel}{n}<\lambda_{t}+\gamma.$

Let
$T_{i}=(D_{f^{i}y}f)^{-1}\circ D_{f^{i}z}f,~0\leq i\leq p-1$,
then
$$D_{z}f^{p}=D_{f^{p-1}y}f\circ T_{p-1}\circ \ldots\circ
D_{fy}f\circ T_1\circ D_{y}f\circ T_{0}.$$

 For $v_{i}\in E_{x}^{i},~1\leq i\leq t$, $x\in \Lambda,$  we define
$\|v_{i}\|_{i}^{'}=\sum_{n=0}^{+\infty}e^{-(\lambda_{i}+2\varepsilon)n}
\|D_{x}f^{n}(v_{i})\|,$ which are clearly  convergent. For
$v=\sum_{i=1}^{t}v_{i},~v_{i}\in E_{x}^{i}$, we define
$\|v\|^{(3)}=\max_{1\leq i\leq t}\{\|v_{i}\|_{i}^{'}\}$. We call
this norm Lyapunov metric number 3. This metric coincides with the
Lyapunov metric number 1 when restricted in the stable bundle, the
direct sum of the sub bundles corresponding to the negative
exponents. This metric is not equivalent to the Riemannain metric in
general. The following estimates are similar to
2.5.\\
$\|Df/_{E_{x}^{i}}\|^{(3)} \leq
e^{\lambda_{i}+2\varepsilon},~1\leq i\leq t;\cdots(3.2.1)$\\
$\frac {1}{\sqrt{d}}\|v\|_{x}\leq \|v\|_{x}^{(3)} \leq
Ce^{\varepsilon k}\|v\|_{x}, ~\forall v\in T_{x}M, ~x\in
\Lambda_{k},\cdots(3.2.2)$\\
 where $C=\frac
{2}{1-e^{-\varepsilon}}$. One can extend this norm to a new norm
 $\|~\|^{(4)}$ by 2.6.

Repeating the process from (3.1.5) to (3.1.11)  in Step 1, we
obtain
$$\parallel D_{f^{i}z}f-D_{f^{i}y}f\parallel^{(4)}\leq
B\eta^{\alpha}e^{-(\alpha_{0}\alpha-\varepsilon)i}, ~0\leq i \leq
p-1,\cdots (3.2.3)$$
$$\parallel D_{f^{i}z}f^{-1}-D_{f^{i}y}f^{-1}
\parallel^{(4)}\leq
B\eta^{\alpha}e^{-(\alpha_{0}\alpha-\varepsilon)(i-1)},~ 1\leq i
\leq p,\cdots (3.2.4)$$ where $B$ is the same constant as in
(3.1.9).
From (3.2.2) and  (3.2.3), for $0\leq i\leq p-1$ we have\\
\begin{eqnarray*}
&\parallel T_i -I
\parallel^{(4)}\\
= &\|(D_{f^{p-i-1}y}f)^{-1}\circ D_{f^{p-i-1}z}f-(D_{f^{p-i-1}y}f)^{-1}
\circ D_{f^{p-i-1}y}f\|^{(4)}\\
\leq &B\eta^{\alpha}e^{-(\alpha_{0}\alpha-\varepsilon)(p-i-1)}
\|(D_{f^{p-i-1}y}f)^{-1}\parallel^{(4)}\\
\leq &CB\eta^{\alpha}e^{-(\alpha_{0}\alpha-2\varepsilon)(p-i-1)}
e^{(k_{0}+2)\varepsilon}\|D_{f^{p-i}y}f^{-1}\parallel.\cdots(3.2.5)\\
\end{eqnarray*}
By using an elementary fact that $\log(1+x)<x,~\forall x>0$, we have\\
 $$\mid\log\parallel
 T_{i}
 \parallel^{(4)}\mid<\parallel T_i-I \parallel ^{(4)}
<CB\eta^{\alpha}e^{-(\alpha_{0}\alpha-2\varepsilon)(p-i-1)}
e^{(k_{0}+2)\varepsilon}\|D_{f^{p-i}y}f^{-1}\parallel.$$
Let $\eta$ be
 small enough such that
$$CB\eta^{\alpha}e^{-(\alpha_{0}\alpha-2\varepsilon)(p-i-1)}
 e^{(k_{0}+2)\varepsilon}\|D_{f^{p-i}y}f^{-1}\parallel
< \frac {\gamma}{5},~0\leq i\leq p-1\cdots (3.2.6).$$

Observe  that for given arbitrarily $ i\in \mathbb{N}$ there
 exists $1\leq j\leq t$ such that $\parallel
D_{f^{i}y}f\parallel^{(3)}=\parallel
D_{f^{i}y}f/_{E_{f^{i}y}^{j}}\parallel^{'}.$
From (3.2.1) and (3.2.6) and the choice of $\eta $ we have\\
\begin{eqnarray*}
&\frac {\log\parallel D_{z}f^{np}\parallel^{(4)}}{pn}
=\frac {\log\parallel (D_{z}f^{p})^{n}\parallel^{(4)}}{pn}
\leq \frac {\log\parallel D_{z}f^{p}\parallel^{(4)}}{p}\\
\leq &\frac {\log(\prod_{i=0}^{p-1}\parallel
D_{f^{i}y}f\parallel^{(3)}\prod_{i=0}^{p-1}\parallel
T_{i}\parallel^{(4)})}{p}\\
= &\frac {1}{p}\log(\prod_{i=0}^{p-1}\parallel
D_{f^{i}y}f\parallel^{(3)})+\frac
{1}{p}\sum_{i=0}^{p-1}\log\parallel
T_{i}\parallel^{(4)}\\
< &\lambda_{t}+2\varepsilon+\frac {\gamma}{5}.
\end{eqnarray*}
Then we have from the choice of $\varepsilon$\\
$$\lim_{n\rightarrow +\infty}\frac
{\log\parallel D_{z}f^{n}\parallel^{(4)}}{n}<\lambda_{t}+\gamma.$$
Note that the norm $\|~\|^{(4)}$
and the Riemannain norm $\|~\|$ are equivalent when restricted on  $Orb(z),$  we get \\
$$\lim_{n\rightarrow +\infty}\frac
{\log\parallel D_{z}f^{n}\parallel}{n}<\lambda_{t}+\gamma.$$

\bigskip

{\bf Step 3} We prove that $\lim_{n\rightarrow +\infty}\frac
{\log\parallel D_{z}f^{n}\parallel}{n}>\lambda_{t}-\gamma.$\\

We now  define another norm, by which we emphasis on the sub bundle
corresponding to the largest Lyapunov exponent $\lambda_t.$ For
$v_{i}\in E_{x}^{i},~1\leq i\leq t-1$, $x\in \Lambda,$ let
$$\|v_{i}\|_{i}^{''}=\sum_{n=0}^{+\infty}
e^{-(\lambda_{i}+2\varepsilon)n} \|D_{x}f^{n}(v_{i})\|;$$ for
$v_{t}\in E_{x}^{t}$, let
$$\|v_{t}\|_{t}^{''}=\sum_{n=0}^{+\infty}e^{(\lambda_{t}-2\varepsilon)n}
\|D_{x}f^{-n}(v_{t})\|.$$ All these  series are clearly convergent.
For $v=\sum_{i=1}^{t}v_{i},~v_{i}\in E_{x}^{i}$, we define
$\|v\|^{(5)}=\max_{1\leq i\leq t}\{\|v_{i}\|_{i}^{''}\}$. We call
this norm Lyapunov metric number 5. The two Lyapunov metrics number
3 and number 5 coincide when restricted on the bundle of the direct
sum of sub bundles corresponding to all but the largest  Lyapunov
exponents. Lyapunov metric number 5 is not equivalent to the
Riemannian metric in general. The following estimates are clear.
\\
$\|Df/_{E_{x}^{i}}\|^{(5)} \leq
e^{\lambda_{i}+2\varepsilon},~1\leq i\leq t-1,
~\|Df/_{E_{x}^{t}}\|^{(5)} \geq e^{\lambda_{t}-2\varepsilon};\quad\cdots(3.3.1)$\\
$\frac {1}{\sqrt{d}}\|v\|_{x}\leq \|v\|_{x}^{(5)} \leq
Ce^{\varepsilon k}\|v\|_{x}, ~\forall v\in T_{x}M, ~x\in
\Lambda_{k},\quad\cdots(3.3.2)$\\
 where $C=\frac
{2}{1-e^{-\varepsilon}}$. One  extends the norm $\|\,\,\|^{(5)}$ to
a new norm  $\|~\|^{(6)}$ by 2.6. Repeating the process from (3.1.5)
to (3.1.11)
in Step 1, we have\\
$$\parallel D_{f^{i}z}f-D_{f^{i}y}f\parallel^{(6)}\leq
B\eta^{\alpha}e^{-(\alpha_{0}\alpha-\varepsilon)i}, ~0\leq i \leq
p-1,\cdots (3.3.3)$$
$$\parallel D_{f^{i}z}f^{-1}-D_{f^{i}y}f^{-1}
\parallel^{(6)}\leq
B\eta^{\alpha}e^{-(\alpha_{0}\alpha-\varepsilon)(i-1)},~ 1\leq i
\leq p,\cdots (3.3.4)$$ where $B$ is  the same constant as in
(3.1.9).

For $\xi>\frac {\frac {2(t-1)\gamma}{5}+1}{1-\frac
{2\gamma}{5}}>1$, let us denote by $K_{\xi}(f^{j}y)$ the following
cones
in $T_{f^{j}y}M,~1\leq j\leq p-1$:\\
$$K_{\xi}(f^{j}y)=\{\sum_{i=1}^{t}v_{i},~v_{i}\in E_{f^{j}y}^{i}, \,\,
1\leq i\leq t; \,\,~\xi\parallel v_{l}
\parallel^{(5)}<\parallel
v_{t}\parallel^{(5)},~1\leq l\leq t-1 \}.$$ From (3.3.1) it follows
that  $D_{f^{j}y}fK_{\xi}(f^{j}y)\subseteq K_{\xi
e^{\lambda_{t}-\lambda_{t-1}-4\varepsilon}}(f^{j+1}y)$.

From (3.3.2), (3.3.3) and (3.3.4),
repeating the proof from (3.2.5) to (3.2.6) in step 2, we have\\
$$\parallel
T_{i}-I
\parallel^{(6)}<\frac {\gamma}{5},~0\leq i\leq p-1\cdots (3.3.5).$$
Now we  consider $T_{j}K_{\xi}(f^{j}y)$. Let $v\in
K_{\xi}(f^{j}y)$, $\,\,v=\sum_{i=1}^{t}v_{i},$ $~v_{i}\in
E_{f^{j}y}^{i},$  $\quad \xi \parallel v_i \parallel ^{(5)}\leq
\parallel v_t\parallel^{(5)},$  $\,\, 1\leq i\leq t-1.$
From (3.3.5) we have\\
$$\frac {\parallel v_{t}\parallel^{(6)}-\sum_{i=1}^{t}\parallel
(T_{j}-I)v_{i}\parallel^{(6)}}{\|v_{t}\|^{(6)}} \geq \frac
{\|v_{t}\|^{(5)}}{\|v_{t}\|^{(5)}}-\sum_{i=1}^{t}\frac
{\gamma}{5}\frac {\parallel v_{i}\parallel^{(5)}}{\|v_{t}\|^{(5)}}
\geq 1-\frac {t\gamma}{5}.\cdots\cdots\cdots \cdots\cdots (3.3.6)$$
And from the choice of $\xi$, for any $1\leq l\leq t-1,$ we have\\
\begin{eqnarray*}
&\parallel v_{t}\parallel^{(6)}-\sum_{i=1}^{t}\parallel
(T_{j}-I)v_{i}\parallel^{(6)}-(\parallel
v_{l}\parallel^{(6)}+\sum_{i=1}^{t}\parallel
(T_{j}-I)v_{i}\parallel^{(6)})\\
\geq &\parallel v_{t}\parallel^{(6)}(1-\frac
{2\gamma}{5}\sum_{i=1}^{t}\frac {\parallel
v_{i}\parallel^{(6)}}{\parallel v_{t}\parallel^{(6)}}-\frac
{1}{\xi})\\
\geq &\parallel v_{t}\parallel^{(6)}(1-\frac {2\gamma}{5}-(\frac
{2(t-1)\gamma}{5}+1)\frac
{1}{\xi})\\
> &0.\\
\end{eqnarray*}
From definition   the norm   $\parallel
\,\,\parallel^{(6)}$ of $T_jv$ coincides with that of the
projection vector $(T_jv)_t$ on $E^t_{f^{j}y},$
$$\parallel T_{j}v\parallel^{(6)}=\parallel
(T_{j}v)_t\parallel_{E_{f^{j}y}^{t}}^{(6)},\,\, \forall\,\, v\in
K_{\xi}(f^{j}y),\,\, \forall \,\, \xi >\frac {\frac
{2(t-1)\gamma}{5}+1}{1-\frac {2\gamma}{5}}>1.\cdots.\cdots\cdots(3.3.7)$$\\
This implies that
$$T_{j}K_{\xi}(f^{j}z)\subseteq
K_1(f^{j}y),\,\,0\leq j\leq p-1.$$
Therefore  we have\\
$$D_{f^{j}y}fT_{j}K_{\xi}(f^{j}z)\subseteq
K_{e^{\lambda_{t}-\lambda_{t-1}-4\varepsilon}}(f^{j+1}y), \,\, 0\leq
j\leq p-1.\cdots(3.3.8)$$
From the choice of $\gamma$ and $\varepsilon$ we have\\
$$e^{\lambda_{t}-\lambda_{t-1}-4\varepsilon}
>\frac
{\frac {2(t-1)\gamma}{5}+1}{1-\frac {2\gamma}{5}}>1, $$ and thus
from (3.3.7) we have
$$\parallel
T_{j+1} v\parallel^{(6)}=\parallel
(T_{j+1}v)_t\parallel_{E_{f^{j+1}z}^{t}}^{(6)},~\forall j\in
\mathbb{N},\,\, v\in
K_{e^{\lambda_{t}-\lambda_{t-1}-4\varepsilon}}(f^{j+1}y), \,\, 0\leq
j\leq p-1.\cdots(3.3.9)$$ From (3.3.6)-(3.3.9)  for any $v\in
K_{\xi}(y)$ it follows
$$\parallel
D_{z}f^{j}v\parallel^{(6)}=\parallel
(D_{z}f^{j}v)_t\parallel_{E_{f^{j}z}^{t}}^{(6)},~\forall j\in
\mathbb{N}$$ and thus by (3.3.1) and (3.3.6) it follows
$$\parallel
D_{z}f^{j}v\parallel^{(6)}\geq (1-\frac {t\gamma}{5})^j
e^{(\lambda_{t}-2\varepsilon)j}\parallel v_t \parallel^{(6)} .$$
Therefore for $v\in K_{\xi}(y)$ we have by the choice of
$\epsilon$
$$\lim_{n\rightarrow +\infty}\frac {1}{n}\log\parallel D_{z}f^{n}v\parallel^{(6)}
\geq \lambda_{t}-2\varepsilon +\log (1-\frac {t\gamma}5)
>\lambda_{t}-(t+1)\gamma.$$
Note that $t\leq d$ we get by replacing $(t+1)\gamma$
by $\gamma$
$$\lim_{n\rightarrow +\infty}\frac {1}{n}\log\parallel D_{z}f^{n}v\parallel^{(6)}
 >\lambda_{t}-\gamma.\cdots\cdots (3.3.10)$$

 Now that  the norm $\|~\|^{(6)}$ and the Riemannian norm $\|~\|$  are equivalent
when  restricted on $Orb(z),$ we complete Step 3 by (3.3.10).
\bigskip

By Step 2 and Step 3 we have
$$\lambda_{t}-\gamma<\lim_{n\rightarrow +\infty}\frac
{\log\parallel D_{z}f^{n}\parallel}{n}<\lambda_{t}+\gamma,$$ by
which we  complete Thmeorem \ref{thm:2}.

\bigskip

\begin{Thm}\label{thm:3}
Let $f:M \rightarrow M$ be a $C^{1+\alpha}$ diffeomorphism of a
compact manifold of dimension $d,$ and let $m$ be an ergodic
hyperbolic measure with Lyapunov exponents
$\lambda_{1}<\cdots<\lambda_{r}<0<\lambda_{r+1}<\cdots<\lambda_{t}(t
\leq d)$. Then the smallest  Lyapunov exponent of $m$ can be
approximated by the smallest  Lyapunov exponents of hyperbolic
periodic orbits. More precisely, for any $\gamma>0$, there exists a
hyperbolic periodic point $z$ with Lyapunov exponents
$\lambda_{1}^{z}\leq\ldots\leq \lambda_{d}^{z}$ such that
$\mid\lambda_{1}-\lambda_{1}^{z}\mid<\gamma$.
\end{Thm}

{\bf Proof} Given $\min_{1\leq i\not=j\leq t}\mid
\lambda_{i}-\lambda_{j}\mid \gg \varepsilon>0$, and for all $k \in
\mathbb{Z}^{+}$, we define
$$\tilde \Lambda_{k}= \tilde
\Lambda_{k}(\{-\lambda_{1},\ldots,-\lambda_{t}\}; \varepsilon)$$
to be all points $x \in M$ for which there is a splitting
$T_{x}M=E_{x}^{1}\oplus \cdot\cdot\cdot\oplus E_{x}^{t}$ with
$$\lim_{n\to \infty}\frac{\log\|df^{n}v\|}{n}=\lambda_{i},~0 \neq
v\in E_{x}^{i}$$ and  with invariant property
$(D_{x}f^{m})E_{x}^{i}=E_{f^{m}x}^{i},~1\leq i\leq t$
and satisfying:\\
\indent$(a)$ $e^{-\varepsilon
k}e^{(-\lambda_{i}-\varepsilon)n}e^{-\varepsilon \mid
m\mid}\leq\|Df^{-n}/_{E_{f^{m}x}^{i}}\| \leq e^{\varepsilon
k}e^{(-\lambda_{i}+\varepsilon)n}e^{\varepsilon \mid m\mid},~1\leq
i\leq t,~\forall m\in
\mathbb{Z},~n\geq 1$;\\

$(b)$ $\tan (Angle(E_{f^{-m}x}^{i},E_{f^{-m}x}^{j})) \geq
e^{-\varepsilon k}e^{-\varepsilon \mid m\mid},~\forall i \neq
j,~\forall m\in
\mathbb{Z}$. \\

We put  $\tilde \Lambda=\tilde
\Lambda(\{-\lambda_{1},\ldots,-\lambda_{t}\};\varepsilon)=\bigcup_{k=1}^{+\infty}
\tilde \Lambda_{k}$ and call $\tilde\Lambda$ a Pesin set. Clearly
$m(\tilde \Lambda)=1.$

 The mesure $m$ is  ergodic and hyperbolic
with respect to $f^{-1},$ for which the Lyapunov exponents are
$$-\lambda_1>...>-\lambda_r>0>-\lambda_{r+1}>...>- \lambda_t.$$
By replacing $f$ by $f^{-1}$ in the proof of Theorem \ref{thm:2}
and by using the Pesin set $\tilde \Lambda$ defined above one can
prove Theorem \ref{thm:3}. We omit the details.

\section {Proof of Theorem \ref{thm:1} }

Based on Theorem \ref{thm:2} and Theorem \ref{thm:3},  we prove
Theorem \ref{thm:1} in this section. We need two more lemmas.

\begin{Lem}\label{lem:4.1}
 Let $f:M
\rightarrow M$ be a $C^{1+\alpha}$ diffeomorphism of a compact
manifold of dimension $d.$ Let $m$ be an ergodic hyperbolic measure
with Lyapunov exponents
$\lambda_{1}<\cdots<\lambda_{r}<0<\lambda_{r+1}<\cdots<\lambda_{t}$
 together with associated splitting
$E^{1}\oplus\cdot\cdot\cdot\oplus E^{t} (t\leq d)$. Then the largest
Lyapunov exponent of $(m,~f^{\wedge^{i}})$, $1\leq i\leq
\sum_{r+1\leq j\leq t}\dim E^{j}$ can be approximated by the largest
Lyapunov exponent of hyperbolic periodic orbits. More precisely, if
we rewrite the Lyapunov spectrum
$\{\lambda_{1},\cdots,\lambda_{t}\}$ of $(m, f)$ as
$\vartheta_{1}\leq \cdots\leq \vartheta_{d}$, then for any
$\gamma>0$, there exists a hyperbolic periodic point $z$ with
Lyapunov exponents $\lambda_{1}^{z}\leq\ldots\leq\lambda_{d}^{z}$
such that
$\mid\sum_{j=d-i+1}^d\vartheta_{j}-\sum_{j=d-i+1}^{d}\lambda_{j}^{z}\mid<\gamma$.
\end{Lem}

{\bf Proof} For all $k \in \mathbb{Z}^{+}$, we define
$$\wedge^{i}_{k}=\wedge^{i}_{k}(\{\sum_{l=1}^{i}\lambda_{j_{l}},~1\leq
j_{1}\leq j_{2}\leq\ldots\leq j_{i}\leq t\};\varepsilon)$$
to be
all points $x \in M$ for which there is a splitting
$$\Lambda^i(x)=\oplus_{1\leq j_{1}\leq j_{2}\leq\ldots\leq
j_{i}\leq t}F_{x}^{j_{1},\ldots,j_{i}},
\,\,\,\,F_{x}^{j_{1},\ldots,j_{i}}=E_{x}^{j_{1}}\wedge\cdots\wedge
E_{x}^{j_{i}}\not=0$$ with
$$\lim_{n\to
\infty}\frac{\log\|Df^{n^{\wedge^{i}}}(v_{j_{1}}\wedge
\ldots\wedge
v_{j_{i}})\|_{\Lambda^{i}}}{n}=\lambda_{j_{1}}+\ldots+\lambda_{j_{i}},
\,\,\,\,v_{j_{1}}\wedge \ldots\wedge v_{j_{i}}\in
F_{x}^{j_{1},\ldots,j_{i}}$$ and with invariant property
$(D_{x}f^{m^{\wedge^{i}}})F_{x}^{j_{1},\ldots,j_{i}}=F_{f^{m}x}^{j_{1},\ldots,j_{i}}$
and satisfying
$$e^{-\varepsilon
k}e^{(\sum_{l=1}^{i}\lambda_{j_{l}}-\varepsilon)n}e^{-\varepsilon
\mid
m\mid}\leq\|Df^{n^{\wedge^{i}}}/_{F_{f^{m}x}^{j_{1},\ldots,j_{i}}}\|_{\wedge^{i}}
\leq e^{\varepsilon
k}e^{(\sum_{l=1}^{i}\lambda_{j_{l}}+\varepsilon)n}e^{\varepsilon
\mid m\mid},$$ $1\leq j_{1}\leq j_{2}\leq\ldots\leq j_{i}\leq
t,~\forall m\in
\mathbb{Z},~n\geq 1.$\\

We put $\Lambda^i=\cup_{k\geq 1}\Lambda^i_k$ and call it a Pesin
set. Clearly $m(\Lambda^i)=1.$ By replacing $f$ by $f^{\Lambda^{i}}$
in the proof of Theorem \ref{thm:2} and by using the Pesin set
$\Lambda^i$ one can prove the  Lemma. We omit the details.

\begin{Lem}\label{lem:4.2}
Let $f:M \rightarrow M$ be a $C^{1+\alpha}$ diffeomorphism of a
compact manifold of dimension $d.$ Let $m$ be an ergodic hyperbolic
measure with Lyapunov exponents
$\lambda_{1}<\ldots<\lambda_{r}<0<\lambda_{r+1}<\ldots<\lambda_{t}$
together with associated splitting $E^{1}\oplus\cdot\cdot\cdot\oplus
E^{t}(t\leq d)$. Then the smallest Lyapunov exponents of
$(m,~f^{\wedge^{i}})$, $1\leq i\leq \sum_{1\leq j\leq r}\dim E^{j}$
can be approximated by the smallest Lyapunov exponent of hyperbolic
periodic orbits. More precisely, if we rewrite the Lyapunov spectrum
$\{\lambda_{1},\cdots,\lambda_{t}\}$ of $(m, f)$ as
$\vartheta_{1}\leq \cdots\leq \vartheta_{d}$, then for $\gamma>0,$
there exists a hyperbolic periodic point $z$ with Lyapunov exponents
$\lambda_{1}^{z}\leq\ldots\leq \lambda_{d}^{z}$  such that
$\mid\sum_{j=1}^{i}\vartheta_{j}-\sum_{j=1}^{i}\lambda_{j}^{z}\mid<\gamma$.
\end{Lem}

{\bf Proof} For all  $k \in \mathbb{Z}^{+},$ we define
$$\tilde \wedge^{i}_{k}=\tilde \wedge^{i}_{k}(\{\sum_{k=1}^{i}(-\lambda_{j_{k}}),~1\leq
j_{1}\leq j_{2}\leq\ldots\leq j_{i}\leq t\};\varepsilon)$$ to be
all points $x \in M$ for which there is a splitting
$$\Lambda^i(x)=\oplus_{1\leq j_{1}\leq
j_{2}\leq\ldots\leq j_{i}\leq t}F_{x}^{'j_{1},\ldots
j_{i}},\,\,\,\,F_{x}^{'j_{1},\ldots
j_{i}}=E_{x}^{j_{1}}\wedge\cdots\wedge E_{x}^{j_{i}}\not=0$$ with
$$\lim_{n\to
\infty}\frac{\log\|df^{-n^{\wedge^{i}}}(v_{j_{1}}\wedge
\ldots\wedge
v_{j_{i}})\|_{\wedge^{i}}}{n}=-(\lambda_{j_{1}}+\ldots+\lambda_{j_{i}}),\,\,\,\,\forall
v_{j_{1}}\wedge \ldots\wedge v_{j_{i}}\in
F_{x}^{'j_{1},\ldots,j_{i}}$$ and with invariant property
$(D_{x}f^{m^{\wedge^{i}}})F_{x}^{'j_{1},\ldots,j_{i}}=
F_{f^{m}x}^{'j_{1},\ldots,j_{i}}$ and satisfying
$$e^{-\varepsilon
k}e^{-(\sum_{l=1}^{i}\lambda_{j_{l}}+\varepsilon)n}e^{-\varepsilon
\mid
m\mid}\leq\|Df^{-n^{\wedge^{i}}}/_{F_{f^{m}x}^{'j_{1},\ldots,j_{i}}}\|_{\wedge^{i}}
\leq e^{\varepsilon
k}e^{-(\sum_{l=1}^{i}\lambda_{j_{l}}-\varepsilon)n}e^{\varepsilon
\mid m\mid},$$ $1\leq j_{1}\leq j_{2}\leq\ldots\leq j_{i}\leq
t,~\forall m\in
\mathbb{Z},~n\geq 1.$\\

We set  $\tilde \Lambda^i=\cup_{k\geq 1}\tilde\Lambda^i_k$ and call
it a Pesin set. Clearly $m(\tilde\Lambda^i)=1.$ By replacing
$f^{-1}$ by $f^{-\Lambda^{i}}$ in the proof of Theorem \ref{thm:3}
and by using the Pesin set $\tilde\Lambda^i$ one can prove Lemma
\ref{lem:4.2}. We omit the details.

\bigskip

{\bf Proof of Theorem \ref{thm:1}} We rewrite the Lyapunov spectrum
$\{\lambda_{1},\cdots,\lambda_{t}\}$ as $\vartheta_{1}\leq
\cdots\leq \vartheta_{d}$. We use   the notations in the proofs of
Theorems\ref{thm:2}- \ref{thm:3} and Lemmas
\ref{lem:4.1}-\ref{lem:4.2} without confusion. For $\forall
\gamma>0$, we can choose $\varepsilon$ following the way in Theorem
\ref{thm:2}. Choose $k_0 \in \mathbb{Z}^{+}$ such that
$$\Gamma_{k_{0}}:=\bigcap_{i=1}^{\sum_{r+1\leq j\leq t}\dim
E^{j}}\wedge^{i}_{k_{0}}\bigcap \Lambda_{k_{0}} \bigcap \tilde
\Lambda_{k_{0}} \bigcap_{i=1}^{\sum_{1\leq j\leq r}\dim
E^{j}}\tilde\wedge^{i}_{k_{0}}$$ has positive $m-$measure. Choose
$\beta(k_{0},\eta,\alpha_{0})>0$ as in Lemma 2.1 and its remark.
According to Theorems \ref{thm:2}- \ref{thm:3} and Lemmas
\ref{lem:4.1}-\ref{lem:4.2} there exists a a hyperbolic point $z\in
\Gamma_{k_{0}}$ with period $p$ and with Lyapunov exponents
$\lambda_{1}^{z}\leq\ldots\leq \lambda_{d}^{z}$
such that  \\
$$\mid\sum_{j=d-i+1}^{d}\vartheta_{j}-\sum_{j=d-i+1}^{d}\lambda_{d}^{z}\mid
<\frac {\gamma}{d},\,\,\,\,~1\leq i \leq \sum_{r+1\leq j\leq
t}\dim E^{i},$$
$$\mid\sum_{j=1}^{i}\vartheta_{j}-\sum_{j=1}^{i}\lambda_{j}^{z}\mid<\frac {\gamma}{d},
\,\,\,\,~1\leq i\leq \sum_{1\leq j\leq r}\dim E^{j}.$$
Thus we have \\
$$\mid\vartheta_{i}-\lambda_{i}^{z}\mid<\gamma,~1\leq i\leq d$$
and complete the proof.
\bigskip

{\bf Acknowledgement} The authors  thank very much Shaobo Gan,
Chao Liang, Geng Liu and Todd Young for their helpful
conversations and the referee for his(her) suggestions.

\section*{ References.}
\begin{enumerate}
\itemsep -2pt

\item [1] L. Barreira, Y. Pesin, Lyapunov exponents and smooth
ergodic theory, Univ. Lect. Ser. 23, AMS, Providence, RI, 2002\\

\item [2] L. Barreira, Y. Pesin, Nonuniform hyperbolicity, dynamics
of systems with nonzero Lyapunov exponents, Cambridge University
Press, Cambridge, 2007\\

\item[3] M. Hirayama, Periodic probability measures are dense in
the set of invariant measures, Dist.  Cont. Dyn. Sys., 9( 2003), 1185-1192\\

\item[4] M. Hirsch, C. Pugh, Stable manifolds and hyperbolic sets,
Proc. Symposia Pure Math. XIV, 133-163, S-S. Chern, S. Smale
Edited, AMS, 1968\\

\item[5] A. Katok, Lyapunov exponents, entropy and periodic orbits
for
diffeomorphisms, Pub. Math. IHES, 51 (1980), 137-173\\

\item[6] A. Katok, L. Mendoza, Dynamical systems with nonuniformly
hyperbolic behavior, Supplement to the book: A. Katok, B.
Hasselblatt, Introduction to the modern theory of
dynamical systems, Cambridge Univ. Press, USA, 1995\\

\item  [7]  C. Liang, G. Liu, W. Sun,  Approximation properties on
invariant measure and Oseledec splitting in non-uniformly
hyperbolic systems, to appear in Trans. Amer. Math. Soc.\\

\item[8] V. I. Oseledec, Multiplicative ergodic theorem, Lyapunov
characteristic numbers for dynamical systems, Trans. Moscow Math.
Soc., 19 (1968), 197-221 \\

\item[9] Y. Pesin, Lyapunov characteristic exponents and ergodic
properties of smooth dynamical systems with an invariant measure,
Sov. Math. Dok., 17(1976), 196-199 \\

\item[10] Y. Pesin, Families of invariant manifolds corresponding
to nonzero Lyapunov exponents, Izvestija, 10(1976), 1261-1305 \\

\item[11] Y. Pesin, Characteristic exponents and smooth ergodic
theory, Russian Mathematical Surveys, 32 no. 4(1997), 55-114\\

 \item[12] M. Pollicott, Lectures on ergodic theory and Pesin
theory on
compact manifolds, Cambridge Univ. Press, 1993 \\

\item[13] D. Ruelle, Ergodic theory of differentiable dynamical
systems, Pub. Math. Lihes., tome 50 (1979), 27-58\\

\item[14] K. Sigmund, Generic properties of invariant measures for
axiom A diffeomorphisms, Invention Math. 11(1970), 99-109

\end{enumerate}
\end{document}